\documentclass[11pt]{article}
\usepackage[margin=1in]{geometry}

\usepackage{setspace}
\doublespace

\usepackage{amsmath}
\usepackage{amssymb}
\usepackage[square]{natbib}
\setcitestyle{numbers}

\newtheorem{theorem}{Theorem}[section]
\newtheorem{corollary}{Corollary}[section]
\newtheorem{lemma}{Lemma}[section]
\newtheorem{remark}{Remark}[section]

\newcommand{\Proof}{\textbf{Proof. }}            
\newcommand{\qed}{$\square$}                     

\newcommand{\RR}{\mathbb{R}}

\newcommand{\ii}{{\mathrm i}}
\newcommand{\ee}{{\mathrm e}}

\newcommand{\dd}{{\mathrm{d}}}

\newcommand{\Pois}{\mathrm{Pois}}    

\newcommand{\Ha}{\mathrm{H}}    
\newcommand{\CP}{\mathrm{CP}} 

\newcommand{\norm}[1]{\|#1\|}                    
\newcommand{\ab}[1]{\vert#1\vert}                
\newcommand{\Ab}[1]{\Big\vert#1\Big\vert}        

\newcommand{\exponent}[1]{\exp\{#1\}}            
\newcommand{\Exponent}[1]{\exp\Bigl\{#1\Bigr\}}  

\newcommand{\eit}{\ee^{\ii t}}                   

\newcommand{\Expect}{\mathrm{E}}                 
\newcommand{\Var}{\mathrm{Var}}                  
\newcommand{\Cov}{\mathrm{Cov}}                  

\newcommand{\pbar}{\overline{p}}

\newcommand{\vfi}{\varphi}
\newcommand{\wE}{\w\Expect}

\newcommand{\eL}{{\cal L}}
\newcommand{\eK}{{\cal K}}

\newcommand{\w}{\widehat}

\newcommand{\binomial}[2]{\genfrac{(}{)}{0pt}{}{#1}{#2}}

\begin{document}
\title{
A Compound Poisson Convergence Theorem for Sums of $m$-Dependent
Variables}

\author{V. \v Cekanavi\v cius  and P. Vellaisamy      \\
{\small
Department of Mathematics and Informatics, Vilnius University,}\\
{\small Naugarduko 24, Vilnius 03225, Lithuania.}\\{\small E-mail:
vydas.cekanavicius@mif.vu.lt } \\{\small and}
\\{\small
 Department of Mathematics, Indian Institute of Technology Bombay,} \\
 {\small Powai, Mumbai-
400076, India.}\\{\small  E-mail: pv@math.iitb.ac.in} }
\date{}

\maketitle

\begin{abstract}

 We prove the Simons-Johnson theorem  for the sums $S_n$ of $m$-dependent random variables,
  with exponential
weights and limiting compound Poisson distribution
$\CP(s,\lambda)$. More precisely, we give sufficient conditions
for $\sum_{k=0}^\infty\ee^{hk}\ab{P(S_n=k)-\CP(s,\lambda)\{k\}}\to
0$ and provide an estimate on the rate of convergence. It is shown
that the Simons-Johnson theorem holds  for weighted Wasserstein
norm as well.
The results are then illustrated for  $N(n;k_1,k_2)$ and
$k$-runs statistics.

\vspace*{.5cm} \noindent {\emph{Key words:} \small   Poisson
distribution, compound Poisson distribution,  m-dependent
variables, Wasserstein norm, rate of convergence.}

\vspace*{.5cm} \noindent {\small {\it MSC 2000 Subject
Classification}:
60F05;  60F15.   
}
\end{abstract}

\newpage

\section{Introduction} 

Simons and Johnson (1971) established an interesting result that
the convergence of the binomial distribution to the limiting
Poisson law can be much stronger than in total variation. Indeed, they proved that if $S_n=X_1+X_2+\cdots+X_n$ has
binomial distribution with parameters $n$, $p=\lambda/n$ and
$g(x)$ satisfies $\sum_0^\infty g(k)\Pois(\lambda)\{k\}<\infty$,
then
\begin{equation}
\sum_{k=0}^\infty g(k)\ab{P(S_n=k)-\Pois(\lambda)\{k\}}\to
0, ~ n \rightarrow \infty,\label{sj71}
\end{equation}
where here and henceforth $\Pois(\lambda)$ denotes Poisson
distribution with mean $\lambda$. The above result was then
extended  to the case of independent and nonidentically
distributed indicator variables by  Chen (1974); see also Barbour
{\it et al.}~(1995) and Borisov and Ruzankin (2002)
 for a comprehensive study in this direction. That similar results
hold for convolutions on measurable Abelian group was proved in
Chen (1975), see also Chen and Roos (1995). Dasgupta (1992) showed
that to some extent,  the binomial distribution in (\ref{sj71})
can be replaced by a negative binomial distribution. Wang (1991)
  later extended Simons and Johnson's result in (\ref{sj71}) to the case of
nonnegative integer valued random variables and compound Poisson limit,
under  the condition
that $P(X_i=k)/P(X_i>0)$ does not depend on $i$ and $n$.

All the above-mentioned works deal with  sums of independent random
variables only. Moreover,  the essential step in the proofs lies in
establishing  an upper bound for the  ratio
$P(S_n=k)/\Pois(\lambda)\{k\}$ or making similar assumptions on
the  measures involved.  The case of dependent random variables is notably less investigated.
  In \v Cekanavi\v cius (2002), the result in (\ref{sj71}) was proved for the
Markov binomial distribution  with $g(k)=\ee^{hk}$. The
possibility to switch from dependent random variables to
independent ones was considered in Ruzankin (2010). However,
results from Ruzankin (2010) are of the intermediate type, since
their estimates usually contain expectations of the unbounded
functionals of the approximated random variables $X_1,\cdots,X_n$,
which still need to be estimated.

In this paper,  we prove the Simons-Johnson theorem with exponential
weights and for the sums of $m$-dependent random variables and
limiting compound Poisson distribution. The main result
contains also estimates on the rate of convergence. A sequence of
random variables $\{X_k \}_{k \geq 1}$
  is called $m$-dependent if, for $1 < s < t < \infty$, $t- s > m$, the sigma-algebras
generated by $X_1,\dots,X_s$ and $X_t, X_{t+1}, \dots$ are
independent. Though the main result is proved for 1-dependent
random variables, it is clear,
 by grouping consecutive summands, that one can reduce the sum of
$m$-dependent variables to the sum of 1-dependent ones. We
exemplify this possibility by considering  $(k_1,k_2)$-events and
$k$-runs.

We consider henceforth the sum $S_n= X_1+X_2+\cdots+X_n$ of
nonidentically distributed 1-dependent random variables
concentrated on nonnegative integers. We denote distribution and
characteristic function of $S_n$ by $F_n(x)$ and $\w F_n(\ii t)$,
respectively.  Note that we include imaginary unit  in the
argument of $\w F_n$, a notation  traditionally preferred over $\w
F_n(t)$ when conjugate distributions  are applied. We define
$j$-th factorial moment of $X_k$ by $\nu_j(k)=\Expect
X_k(X_k-1)\cdots(X_k-j+1)$, $k=1,2,\dots,n$; $j=1,2,\dots$~.  Let
\begin{equation*} \Gamma_1=\Expect
S_n=\sum_{k=1}^n\nu_1(k),\quad \Gamma_2=\frac{1}{2}(\Var
S_n-\Expect
S_n)=\frac{1}{2}\sum_{k=1}^n\big(\nu_2(k)-\nu_1^2(k)\big)+\sum_{k=2}^n\Cov(X_{k-1},X_k).
\label{gamma}
\end{equation*}
Formally,
\begin{equation}\label{formal}
 \w F_n(\ii t)=\exponent{\Gamma_1(\eit-1)+\Gamma_2(\eit-1)^2+\dots}.
\end{equation}
It is clear that  Poisson limit occurs  only if
$\Gamma_1\to\lambda$, $\Gamma_2\to 0$, and other factorial
cumulants also tend to zero. Similar arguments apply for compound
Poisson limit as well.

 Next, we introduce
compound Poisson distribution
$\CP(s,\lambda)=\CP(s,\lambda_1,\dots,\lambda_s)$, where $s\geqslant 1$ is an integer.
Let $N_i$ be independent Poisson random variables with parameters
$\lambda_i\geqslant 0$, $i=1,2,\dots,s$. Then $\CP(s,\lambda)$ is defined as the
 distribution of $N_1+2N_2+3N_3+\dots+sN_s$ with characteristic function
\begin{equation}
\w\CP(s,\lambda)(\ii t)=\Exponent{\sum_{m=1}^s\lambda_m(\ee^{\ii
tm}-1)}=\Exponent{\sum_{j=1}^s(\eit-1)^j\sum_{m=j}^s\binomial{m}{j}\lambda_m}.
\label{CPch}
\end{equation}

\noindent Note also that
\begin{equation*}
N_1+2N_2+\dots+sN_s \stackrel{{\cal L}}{=} Y_1+ Y_2+\dots+Y_N,
\end{equation*}
where the $Y_j$ are independent random variables with $P(Y_1=j)=
\lambda_j/(\sum_{i=1}^{s} \lambda_i),$ for $1 \leq j \leq s$ and
$N \sim \Pois(\sum_{i=1}^{s} \lambda_i).$ It is clear that when
$s=1$, $\CP(1,\lambda)=\Pois(\lambda)$, the distribution of $N$ in
this case.

Let $M$ be a signed measure concentrated on nonnegative integers.
The total variation norm of $M$ is denoted by
$\norm{M}=\sum_{m=0}^\infty\ab{M\{m\}}$. Properties of the  norm
are discussed in detail in Shiryaev (1995), pp. 359--362. The
total variation norm is arguably the most popular metric used for
estimation of the accuracy of approximation of discrete random
variables. The Wasserstein (or Kantorovich) norm is defined as
$\norm{M}_W= \sum_{m=0}^\infty\Ab{\sum_{k=0}^mM\{k\}}$.  For other
expressions of $\norm{M}$ and $\norm{M}_W$ one can consult
appendix A1 in Barbour et al. (1992).

\section{The Main Results}

Henceforth, we assume that all random variables are uniformly
bounded from above, that is, $X_i\leqslant C_0, 1 \leqslant i
\leqslant n$. Here, $C_0\geqslant 1$ is some absolute constant.
First, we formulate sufficient conditions for
 compound Poisson limit with exponential
weights.

\begin{theorem}\label{converg} Let  $X_i$ be nonidentically distributed 1-dependent random variables
concentrated on nonnegative integers, $X_i\leqslant C_0, 1
\leqslant i \leqslant n$. Let $F_n(x)$ denote the distribution of
$S_n= X_1+X_2+\cdots+X_n$ and let $\CP(s,\lambda)$ be defined by
(\ref{CPch}). Let $s \geqslant 1$ be an integer, $\lambda_j
\geqslant 0, 1 \leqslant j \leqslant s$, and $h\geqslant 0$ be
fixed numbers. If, as $n\to\infty$,
\begin{eqnarray}
&&\max_{1\leqslant j\leqslant n}\nu_1(j)\to 0,\label{s1}\\
&&\frac{1}{m!}\sum_{j=1}^n\nu_m(j)\to\sum_{l=m}^s\binomial{l}{m}\lambda_l,\quad
m=1,2,\dots,s;\label{s2}\\
&&\sum_{j=1}^n\nu_{s+1}(j)\to 0,\label{s3}\\
&&\sum_{j=2}^n\ab{\Cov(X_{j-1},X_j)}\to 0,\label{s4}
\end{eqnarray}
then
 \begin{equation}
 \sum_{k=0}^\infty\ee^{hk}\ab{F_n\{k\}-\CP(s,\lambda)\{k\}}\to 0.
\label{conv}
\end{equation}
\end{theorem}
\begin{remark} (i) Assumption $C_0\geqslant 1$ is not restrictive.
Indeed, $X_i<1$ is equivalent to the trivial case $X_i\equiv 0$,
since we assume that $X_i$ is concentrated on integers.

(ii) Technical assumption that all random variables are uniformly
bounded significantly simplifies all proofs. Probably it can be
replaced by some  more general uniform smallness conditions for
the tails of distributions.

(iii) Conditions for convergence to compound Poisson distribution
can be formulated in various terms. In Theorem \ref{converg} we
used factorial cumulants. Observe that such approach allows
natural comparison of the characteristic functions due to the
exponential structure of  $\w\CP(s,\lambda)(\ii t)$.

(iv) Assumptions (\ref{s1})--(\ref{s4}) are sufficient for
convergence, but not necessary.
For example, consider the case $s=2$ and
compare (\ref{formal}) and (\ref{CPch}). The convergence then
implies $\Gamma_1\to\lambda_1+2\lambda_2$ and
$\Gamma_2\to\lambda_2$. If we assume, in addition (\ref{s1}), then
the last condition is equivalent to
\[\frac{1}{2}\sum_{j=1}^n\nu_2(j)+\sum_{j=2}^n\Cov(
X_{j-1},X_j)\to\lambda_2,\] and is more general than the
assumptions $\sum_1^n\nu_2(j)/2\to \lambda_2$ and (\ref{s4}).
\end{remark}

Observe that we can treat (\ref{sj71}) as a weighted total variation norm with
increasing weights. A natural question that arises is the following: is it possible to
extend this result to stronger norms? If we consider the
Wasserstein norm, then the answer is affirmative, see Lemma
\ref{WSJ} below. Let $F_n(k)=F_n\{[0,k]\}$ and
$\CP(s,\lambda)(k)=\CP(s,\lambda)\{[0,k]\}$ denote the corresponding
distribution functions. For exponentially weighted Wasserstein
norm, we have the following inequality:
\begin{equation}
\sum_{k=0}^\infty\ee^{hk}\ab{F_n(k)-\CP(s,\lambda)(k)}
\leqslant\frac{1}{\ee^h-1}
\sum_{k=0}^\infty\ee^{hk}\ab{F_n\{k\}-\CP(s,\lambda)\{k\}},
\label{ws}
\end{equation}
provided the left-hand side is finite and $h>0$. We see that,
though Wasserstein norm (which corresponds to the case $h=0$) is
stronger than the total variation norm, the weighted Wasserstein
norm is bounded from above by the correspondingly weighted total variation norm.
Consequently, from (\ref{ws}) and Theorem \ref{converg}, the
following corollary immediately follows.

\begin{corollary}\label{wteorema} Let
$\lambda_1\geqslant 0,\dots,\lambda_s\geqslant 0,$ and $
s\geqslant 1$ be an integer. Assume   conditions
(\ref{s1})--(\ref{s4}) are satisfied. Then, for fixed $h> 0$,
\begin{equation}
\sum_{k=0}^\infty\ee^{hk}\ab{F_n(k)-\CP(s,\lambda)(k)}\to 0.
\end{equation}
\end{corollary}

\noindent  Indeed, Theorem \ref{converg} follows from more general Theorem \ref{teorema} given below. Assuming
$\max_j\nu_1(j)$ to be small, but not necessarily  converging to zero,
we obtain estimates of remainder terms. Let
\begin{eqnarray} \label{eqn11}
a&=&a(h,C_0)=\ee^{hC_0}(2+h)\sqrt{C_0},\quad \psi=
\Exponent{\max\Bigl(4a^2\Gamma_1,\sum_{m=1}^s\lambda_m(\ee^{hm}+1)\Bigr)},\\ \nonumber
K_1&=&\psi\sqrt{\pi+1}(\ee^h+1)^{s}(s+1+4a^2\Gamma_1),\quad
K_2=\psi\sqrt{\pi+1}(s+1+4a^2\Gamma_1)\frac{\ee^{hC_0}(\ee^h+1)^{s+1}}{(s+1)!},\\ \nonumber
K_3&=&16\psi a^4\sqrt{\pi+1}(5+6a^2\Gamma_1),\quad
K_4=4\psi a^3\sqrt{\pi+1}(1.1+a^2\Gamma_1).
\end{eqnarray}

\noindent Let us denote henceforth $ \nu_1^{(n)}= \max_{1\leqslant j\leqslant n}\nu_1(j),$ for simplicity.
We are ready to state the main result of this paper.
\begin{theorem}\label{teorema} Let $s \geqslant 1$ be an integer, $h\geqslant 0$, $\lambda_j \geqslant 0, 1 \leqslant j \leqslant s$, and
 let $a^2 \nu_1^{(n)} \leqslant
1/100$. Then,
\begin{eqnarray} \label{neqn12}
\lefteqn{\sum_{k=0}^\infty\ee^{hk}\ab{F_n\{k\}-\CP(s,\lambda)\{k\}}\leqslant
K_1\sum_{m=1}^s\Ab{\frac{1}{m!}\sum_{j=1}^n\nu_m(j)-\sum_{l=m}^s\binomial{l}{m}\lambda_l}}\hskip 3cm\nonumber\\
&&+K_2\sum_{j=1}^n\nu_{s+1}(j)+K_3\sum_{j=1}^n\nu_1^2(j)+K_4\sum_{j=2}^n\ab{\Cov(X_{j-1},X_j)}.
\label{main}
\end{eqnarray}
\end{theorem}
We next illustrate the results for the cases  $s=1$ and $s=2$, which are of particular interest. Note here
 the corresponding  limiting distributions are as follows:

 \vspace*{-0.3cm}
\[\w\Pois(\lambda)(\ii t)=\exponent{\lambda(\eit-1)},\quad
\w\CP(2,\lambda)(\ii
t)=\exponent{\lambda_1(\eit-1)+\lambda_2(\ee^{2\ii t}-1)}.
\]

\vspace*{-0.3cm}
The following corollary is immediate from (\ref{neqn12}).

\begin{corollary} \label{cor2.2} Let $a^2 \nu_1^{(n)} \leqslant
1/100$. Assume $h\geqslant 0$,  $\lambda,\lambda_1$ and
$\lambda_2$ are positive reals. Then,
\begin{eqnarray}
 (i)~~ \lefteqn{\sum_{k=0}^\infty\ee^{hk}\ab{F_n\{k\}-\Pois(\lambda)\{k\}}}\hskip
0.5cm\nonumber\\
&\leqslant&
C_1(h,\lambda)\exponent{4a^2\Gamma_1}\bigg\{\ab{\Gamma_1-\lambda}+\sum_{j=1}^n\nu_2(j)
+\sum_{j=1}^n\nu_1^2(j)+\sum_{j=2}^n\ab{\Cov(X_{j-1},X_j)}\bigg\},
\label{poison}\\
(ii)~~ \lefteqn{\sum_{k=1}^n\ee^{hk}\ab{F_n\{k\}-\CP(2,\lambda)\{k\}}}\hskip
0.5cm\nonumber\\
&\leqslant&
C_2(h,\lambda_1,\lambda_2)\exponent{4a^2\Gamma_1}\bigg\{
\ab{\Gamma_1-\lambda_1-2\lambda_2}+\Ab{\sum_{j=1}^n\nu_2(j)-2\lambda_2}
+\sum_{j=1}^n\nu_3(j)\nonumber\\
&&+\sum_{j=1}^n\nu_1^2(j)+\sum_{j=2}^n\ab{\Cov(X_{j-1},X_j)}\bigg\}.\label{CP2}
\end{eqnarray}
Note here the constants  $C_1$ and  $C_2$ depend  on
 $h$, $\lambda,\lambda_1$ and $\lambda_2$ only.
\end{corollary}
\begin{remark} (i) Applying (\ref{ws}), we can obtain the
estimate for exponentially weighted Wasserstein norm, similar to Theorem \ref{teorema}.

(ii) Let us consider the sum of independent Bernoulli variables,
$W=\xi_1+\dots+\xi_n$, where $P(\xi_i=1)=1-P(\xi_i=0)=p_i$. Assume
 that, for some fixed $\lambda>0$, the parameter $p_i$ satisfies $\sum_{k=1}^np_i=\lambda$
and $\sum_{k=1}^np_i^2\to 0$, as $n\rightarrow \infty.$ Then, putting $h=0$ in
(\ref{poison}), we obtain an estimate for total variation metric as
\[\sum_{k=0}^\infty\ab{P(W=k)-\Pois(\lambda)\{k\}}\leqslant
C_3\sum_{j=1}^np_j^2,\] if $n$ is sufficiently large. Observe that this
estimate is of the right order.
\end{remark}

\noindent We next show that Simons-Johnson result holds for
convergence associated  with $(k_1, k_2)$-events and $k$-runs,
which have applications in statistics. For example, the number of
$k$-runs have been used to develop certain nonparametric tests for
randomness. See O'Brien and Dyck (1985) for more details.

%
\section{Some Examples}
In examples below, we assume $\lambda,\lambda_1,\lambda_2$ and $h\geqslant 0$ are
some absolute constants.

\noindent\textbf{ 1. Number of $(k_1,k_2)$ events.} Consider a
sequence of independent Bernoulli trials with the same success
probability $p$. We say that $(k_1,k_2)$-event has occurred if
$k_1$ consecutive failures are followed by $k_2$ consecutive
successes. Such sequences can be meaningful in biology (see Huang
and Tsai (1991), p. 126), or in agriculture, since sequences of
rainy and dry days have impact on the yield of raisins (see Dafnis
et al. (2010), p. 1698).

More formally, let $\eta_i$ be independent Bernoulli
 $Be(p)$ ($0<p<1$)  variables and
$Z_j=(1-\eta_{j-m+1})\cdots
(1-\eta_{j-k_2})\eta_{j-k_2+1}\cdots\eta_{j-1}\eta_j$,
$j=m,m+1,\dots,n$, where $m=k_1+k_2$ and $k_1>0$ and $k_2>0$ are
fixed integers. Then, $N(n;k_1,k_2)=Z_m+Z_{m+1}+\cdots+Z_n$
denotes the number of $(k_1, k_2)$ events in
 $n$ Bernoulli trials. We denote the distribution
of $N(n;k_1,k_2)$ by $\Ha$.  It is well known that $N(n;k_1,k_2)$
has limiting Poisson distribution, see Huang and Tsai (1991) and
Vellaisamy (2004). Note also that $Z_1,Z_2,\dots$ are
$m$-dependent. Consequently, the results of previous section
cannot be applied directly. However, one can group the summands in
the following natural way: \vspace*{-0.3cm}
\begin{eqnarray*}
N(n;k_1,k_2)&=&(Z_m+Z_{m+1}+\cdots+Z_{2m-1})+(Z_{2m}+Z_{2m+1}+\cdots+Z_{3m-1})+\dots \\
             &=&X_1+X_2+\dots~.
\end{eqnarray*}

\vspace*{-0.3cm}
\noindent Here, each $X_j$, with probable exception of the last one, contains $m$
 summands. Let $K$ and $\delta$ be the
integer and fractional parts of $(n-m+1)/m$, respectively, so that
\begin{equation*}
K=\bigg\lfloor \frac{n-m+1}{m}\bigg\rfloor,\quad
\frac{n-m+1}{m}=K+\delta,\quad 0\leqslant\delta<1, \label{Kdelta}
\end{equation*}
and $a(p)=(1-p)^{k_1}p^{k_2}$. Then, considering the structure of
new variables $X_j$ we see that, for $j=1,\dots, K$
\begin{equation*}
X_j= \begin{cases} 1, & \mbox{with probability } ma(p), \\
0, & \mbox{with probability } 1-ma(p), \end{cases}\quad
X_{K+1}= \begin{cases} 1, & \mbox{with probability } \delta ma(p), \\
0, & \mbox{with probability } 1-\delta ma(p). \end{cases}
\end{equation*}
 Consequently, $\nu_2(j)=\nu_2(K+1)=0$,
$\nu_1(j)=ma(p)$, $\nu_1(K+1)=\delta a(p)$, $\Gamma_1=(n-m+1)a(p)$
and  we obtain, checking for nonzero products,
 \[\Expect
 (X_1X_2)=a^2(p)(m+(m-1)+(m-2)+\dots+1)=\frac{a(p)^2m(m+1)}{2},\quad \Expect (X_{K}X_{K+1})=\frac{\delta m(\delta m+1)a^2(p)}{2}.\]
Therefore,
\[\Cov(X_{j-1},X_j)=-\frac{m(m-1)a^2(p)}{2},\quad
\Cov(X_{K},X_{K+1})=\frac{a^2(p)\delta m(\delta m+1-2m)}{2},\]
for $j=1,2,\dots, K$.
 Consequently, if $(n-m+1)a(p)\to\lambda$, then
\begin{equation*}
\sum_{j=0}^\infty\ee^{hj}\ab{\Ha\{j\}-\Pois(\lambda)\{j\}}\to
0.\label{HaP}
\end{equation*}
Indeed, we have $a(p)=o(1)$ and
\begin{eqnarray*}\sum_{j=2}^{K+1}\ab{\Cov(X_{j-1},X_j)}&\leqslant&
\frac{Km(m-1)a^2(p)+a^2(p)\delta m(2m-1-\delta m)}{2}\\
&\leqslant& (Km+\delta m)a^2(p)m=(n-m+1)a^2(p)\to 0,\\
\sum_{j=1}^{K+1}\nu_1^2(j)&\leqslant& a(p)\Gamma_1\to 0
\end{eqnarray*}
Using (\ref{poison}) of Corollary \ref{cor2.2}, we see that (\ref{conv}) holds with
$\CP(1,\lambda)$.

\vspace*{0.4cm} \noindent\textbf{ 2. Statistic of $k$-runs.} Let
$\eta_i, 1 \leq i \leq n+k-1, $ be independent Bernoulli $Be(p)$
($0<p<1$) variables and let
$Z_j=\eta_j\eta_{j+1}\cdots\eta_{j+k-1}$. Then
$S=Z_1+Z_2+\dots+Z_n$ is called $k$-runs statistic. Runs
statistics are important in reliability theory ($m$ consecutive
$k$ out of $n$ failure system) and quality control (see, for
discussion, Wang and Xia (2008)). Approximations of 2 or  $k$-runs
statistic (including the case of different probabilities $p_i$) by
various distributions have been considered in numerous papers, see
R\"{o}llin (2005)  and Wang and Xia (2008) and the references
therein. As in the previous example, we switch from $k$-dependent
case to 1-dependent one by grouping $k$ consecutive summands as
$X_1=Z_1+\dots+Z_k$, $X_2=Z_{k+1}+\cdots+Z_{2k}$ and so on. Note
that such a grouping is not unique. For example, it is possible to
group ($k$-$1$) consecutive summands. Let $K$ denote the integer
part of $(n/k)$, where $k$ is fixed. Next, we apply Corollary
\ref{cor2.2}. It is obvious that $\Gamma_1=np^k$,
$\nu_2(K+1)=o(1)$, and $E(X_{K}X_{K+1})=o(1)$ as $n\to\infty$. For
$j=2,\dots,K$, we have $\Expect (X_{j-1}X_j)\leqslant C(k)p^{k+1}$
and $\nu_2(j)\leqslant C(k)p^{k+1}$.
 Indeed, in both the cases, at
least two of $Z_i$'s must be equal to unity. Next, note that

\vspace*{-0.4cm}
\begin{equation*}
\sum_{j=2}^K\ab{\Cov(X_{j-1},X_j)}\leqslant\sum_{j=2}^K\Expect
(X_{j-1}X_j)+\sum_{j=2}^K\nu_1(j-1)\nu_1(j)\leqslant C(k)np^{k+1}.
\end{equation*}
 Consequently, if $np^k\to\lambda$, then  (\ref{conv}) holds for
$F_n=\eL(S)$ with limiting $ Pois(\lambda)$ distribution.

\noindent \textbf{3. Convergence to $\CP(2,\lambda)$.} By slightly
modifying 2-runs, we construct an example of 1-dependent summands with
limiting compound Poisson distribution. Let  $\eta_i\sim
 Be(p)$,  ($0<p<1$, $i=1,\dots,n+1$) and  $\xi_j\sim Be(\pbar)$, ($0<\pbar<1$, $j=1,\dots,n$) be
  two sequences of independent
 Bernoulli variables (any $\xi_j$ and $\eta_i$ are also independent). Let $X_1=\eta_1\eta_2+2\xi_1(1-\eta_1\eta_2)$,
 $X_2=\eta_2\eta_3+2\xi_2(1-\eta_2\eta_3)$, $X_3=\eta_3\eta_4+2\xi_3(1-\eta_3\eta_4)$ and so on. Let $S=X_1+\cdots+X_n$. It is
 obvious  that $X_1,X_2,\dots, X_n$ are 1-dependent random
 variables.  Moreover,
 \begin{eqnarray*}
 \nu_1(j)&=&p^2+2\pbar(1-p^2),\quad \nu_2(j)=\Expect
 (X_j(X_j-1))=2\pbar(1-p^2),\\
 \nu_3(j)&=&0,\quad\ab{\Cov(X_1,X_2)}\leqslant\Expect
 (X_1X_2)+\nu_1(1)\nu_1(2)\leqslant C(p^3+p^2\pbar+\pbar^2).
 \end{eqnarray*}
Let $ np^2\to\lambda_1$ and $n\pbar\to\lambda_2$, as $n \rightarrow \infty.$ Then
\[\Gamma_1=n\nu_1(1)\to\lambda_1+2\lambda_2,\quad
\sum_{j=1}^n\nu_2(j)\to 2\lambda_2,\quad\sum_{j=1}^n\nu_1^2(j)\to
0, \quad\sum_{j=2}^n\ab{\Cov(X_{j-1},X_j)}\to 0.\]
Therefore,  it follows from
(\ref{CP2}) that
\[
\sum_{k=0}^n\ee^{kh}\ab{P(S=k)-\CP(2,\lambda)\{k\}}\to 0,
\]
leading to strong compound Poisson convergence.

\section{Proofs}

It is not difficult to observe that the weighted sum in Theorem
\ref{teorema}  can be treated as the total variation of some
conjugate measure. Indeed,
\[
\sum_{k=0}^\infty\ee^{hk}\ab{F_n\{k\}-\CP(s,\lambda)\{k\}}=\norm{M}.\]
Here, $ M\{k\}=\ee^{hk}(F_n\{k\}-\CP(s,\lambda)\{k\})$. For
estimation of $\norm{M}$ we  apply the characteristic function
method. Observe that  $\w M(\ii t)=\w F_n(\ii
t+h)-\w\CP(s,\lambda)(\ii t+h)$. We need to estimate $\ab{\w M(\ii
t)}$. Therefore, the crucial step in the proof, is expansion of
$\w F_n(\ii t+h)$ in moments of $S_n$. The essential tool for this
 is Heinrich's (1982) representation of $\w F_n(\ii t)$ as a
product of $n$ functions. For Heinrich's representation, we need
some additional notations.

\noindent  Let $\{U_{k}\}_{k \geq 1}$ be a sequence of arbitrary
real or complex-valued random variables. Also, let  $\w\Expect
(U_1)=\Expect (U_1)$ and, for $k\geqslant 2$, define
\begin{equation*}
 \w\Expect (U_1,U_2,\cdots, U_k)=\Expect (U_1U_2\cdots
U_k)-\sum_{j=1}^{k-1}\w\Expect (U_1,\cdots ,U_j)\Expect
(U_{j+1}\cdots U_{k}). \label{capZ}
\end{equation*}
Then, it is obvious that $\wE(X_{k-1},X_k)=\Cov(X_{k-1},X_k)$.


We require the following two lemmas from Heinrich (1982).

\begin{lemma}(Heinrich (1982)) \label{Hei3aa} Let $U_1,U_2,\dots,U_k$ be 1-dependent complex-valued random variables with
$\Expect\ab{U_m}^2<\infty$, $ 1 \leq m \leq k. $ Then,
\[
\ab{\wE (U_1, U_2, \cdots, U_k)}\leqslant
2^{k-1}\prod_{m=1}^k(\Expect\ab{U_m}^2)^{1/2}.
\]
\end{lemma}
For a complex number $z$, let
\[
w(z)=\max_{1\leqslant k\leqslant
n}\sqrt{\Expect\ab{\ee^{zX_k}-1}^2},\quad \eK=\{z: w(z)\leqslant
1/6\}.
\]

\begin{lemma}\label{Horiginal} (Heinrich (1982)) Let
$X_1,X_2,\dots,X_n$ be a sequence of 1-dependent random variables.
Then for each $z\in\eK$, the following product representation holds:
\begin{equation*}
\Expect (\ee^{zS_n})=\vfi_1(z)\vfi_2(z)\cdots\vfi_n(z).
\label{product}
\end{equation*}
Here, $\vfi_1(z)=\Expect (\ee^{zX_1})$ and for $k=2,\dots,n$,
\begin{equation}\vfi_k(z)= 1+ \Expect
(\ee^{zX_k}-1)+\sum_{j=1}^{k-1}\frac{\wE \Big((\ee^{zX_j}-1),(\ee^{zX_{j+1}}-1),\dots,
(\ee^{zX_k}-1)\Big)}{\vfi_j(z)\vfi_{j+1}(z)\dots
\vfi_{k-1}(z)},\label{repr}\end{equation}
Further,
\begin{equation}
\ab{\vfi_k(z)-1}\leqslant\ab{\Expect(\ee^{zX_k})-1}+\frac{2 \Big(\Expect\ab{\ee^{zX_{k-1}}-1}^2\Expect
\ab{\ee^{zX_k}-1}^2\Big)^{1/2}}{1-4w(z)}. \label{fir}, \end{equation}
\end{lemma}
for $z\in\eK$ and $ 1 \leqslant k \leqslant n.$

\vspace*{0.4cm}
\noindent In addition, we use the following notation: $u=\ii t+h$,
$Y_j=\exponent{uX_j}-1$, $\Psi_{jk}=\wE(Y_j,\dots,Y_k)$. We use
symbol $\theta$ to denote a real or a complex number
satisfying $\ab{\theta}\leqslant 1$.
 Assume
 $\nu_j(k)=0$ and $X_k=0$ for $k\leqslant 0$. Moreover,
 $\nu_j(k)=0$ if $X_j<k$ and $\nu_1^{(n)}=\max_{1\leqslant
 j\leqslant n}\nu_1(j)$, as earlier.
The primes denote the derivatives with respect to $t$.

\begin{lemma} The following relations hold for all $t$, $k=1,\dots,n$, and an integer $s\geqslant 1$:
\begin{eqnarray}
\ab{Y_k}&\leqslant&\ee^{hC_0}(2+h)X_k,\quad\ab{Y_k}^2\leqslant
a^2X_k,\quad\Expect\ab{Y_k}\leqslant
a\nu_1(k),\quad\Expect\ab{Y_k}^2\leqslant a^2\nu_1(k),\label{Yj}\\
\ab{Y_k'}&\leqslant&\ee^{hC_0}X_k,\quad\ab{Y_k'}^2\leqslant\ee^{2hC_0}C_0X_k,
\quad\Expect\ab{Y_k'}\leqslant\frac{a}{2}\nu_1(k),\quad\Expect\ab{Y_k'}^2\leqslant\frac{a^2}{4}\nu_1(k),
\label{Yj0}\\
 \Expect Y_k &=&
\sum_{m=1}^s\frac{\nu_m(k)}{m!}(\ee^u-1)^m+\theta\ee^{hC_0}(\ee^h+1)^{s+1}\frac{\nu_{s+1}(k)}{(s+1)!},\label{Yj1}\\
\Expect
Y_k'&=&\ii\sum_{m=1}^s\frac{\nu_m(k)}{(m-1)!}\ee^u(\ee^u-1)^{m-1}+\theta\ee^{hC_0}(\ee^h+1)^s\frac{\nu_{s+1}(k)}{s!}.\label{Yj2}
\end{eqnarray}
\end{lemma}

\noindent \Proof Since $\ab{\exponent{\ii t (X_k-j)}}=1$,   we have
\begin{eqnarray*}\ab{Y_k}&\leqslant&\ee^{hX_k}\ab{\ee^{\ii t
X_k}-1}+\ee^{hX_k}-1\leqslant \ee^{hX_k} \Big(\ab{\ee^{\ii
t(X_k-1)}}+\ab{\ee^{\ii
t(X_k-2)}}+\cdots+1 \Big)\ab{\eit-1}\\
&&+hX_k\ee^{hX_k}\leqslant
\ee^{hC_0}X_k(\ab{\eit}+1)+hX_k\ee^{hC_0}\\
&\leqslant &
\ee^{hC_0}(2+h)X_k.
\end{eqnarray*}
Other relations of (\ref{Yj}) now follow. The proof of (\ref{Yj0})
is obvious. For the proof of (\ref{Yj1}), we apply Bergstr\"om
(1951) identity
\begin{equation}
\alpha^N\,=\,\sum_{m=0}^s\binomial{N}{m}\beta^{N-m}(\alpha-\beta)^m
+\sum_{m=s+1}^N\binomial{m-1}{s}\alpha^{N-m}(\alpha-\beta)^{s+1}\beta^{m-s-1},
            \label{bergident}
\end{equation}
which holds for any  numbers $\alpha,\beta$ and $s=0,1,2,\dots,N$. Let
 $\binomial{j}{k}=0$, for $k>j$. Then, (\ref{bergident})
holds for all $s=0,1,\dots$ .  We apply (\ref{bergident}) with
$N=X_k$, $\alpha=\ee^{u}$ and  $\beta=1$. Then,
\begin{equation}
Y_k\,=\,\sum_{m=1}^s\binomial{X_k}{m}(\ee^{u}-1)^m
+\sum_{m=s+1}^{X_k}\binomial{m-1}{s}\ee^{u(X_k-m)}(\ee^u-1)^{s+1}.
            \label{berg}
            \end{equation}
Using the results
\[\sum_{m=s+1}^N\binomial{m-1}{s}=\binomial{N}{s+1},\quad\ab{\ee^u}=\ee^h,\]
we obtain
\[\sum_{m=s+1}^{X_k}\binomial{m-1}{s}\ab{\ee^{u(X_j-m)}}\leqslant\ee^{hC_0}\binomial{X_j}{s+1}.
\]
The proof of (\ref{Yj1}) now follows by finding the mean of $Y_k$ in (\ref{berg}) and using the definition
of $\nu_j(k)$.

\noindent For the proof of (\ref{Yj2}),
we once again apply (\ref{bergident}) to obtain
 \begin{eqnarray*}
Y_k'&=&\ii X_k\ee^{uX_k}=\ii X_k\ee^{u}\ee^{u(X_k-1)}\\
&=& \ii
X_k\ee^u\bigg\{\sum_{m=0}^{s-1}\binomial{X_k-1}{m}(\ee^u-1)^m+(\ee^u-1)^s\sum_{m=s}^{X_k-1}\binomial{m-1}{s-1}
\ee^{u(X_k-1-m)}\bigg\}.
\end{eqnarray*}
The rest of the proof is the same as that of  (\ref{Yj1}) and,
therefore, omitted. \qed

\begin{lemma} \label{Psi} Let $a^2 \nu_1^{(n)} \leqslant 0.01$. Then, for $k=4,\dots,n$ and $j=1,\dots,k-3$,
\[
\ab{\Psi_{jk}}\leqslant 250
a^4\bigg(\frac{1}{5}\bigg)^{k-j}\sum_{l=0}^3\nu_1^2(k-l),\quad
\ab{\Psi_{jk}'}\leqslant 125a^4(k-j+1)
\bigg(\frac{1}{5}\bigg)^{k-j}\sum_{l=0}^3\nu_1^2(k-l)\]
   and for $k=2,\dots,n$; $j=1,\dots,k-1$,
\begin{eqnarray*}
   \ab{\Psi_{jk}} & \leqslant &
   5a^2\bigg(\frac{1}{5}\bigg)^{k-j}[\nu_1(k-1)+\nu_1(k)], \\
   \ab{\Psi_{jk}'} & \leqslant &
   (2.5)a^2(k-j+1)\bigg(\frac{1}{5}\bigg)^{k-j}[\nu_1(k-1)+\nu_1(k)].
\end{eqnarray*}
\end{lemma}
\Proof From Lemma \ref{Hei3aa} and (\ref{Yj}), we have
\[\ab{\Psi_{jk}}\leqslant
2^{k-j}\prod_{l=j}^k\sqrt{a^2\nu_1(l)}\leqslant
2^{k-j}(0.1)^{k-j-3}a^4\sqrt{\nu_1(k)\nu_1(k-1)\nu_1(k-2)\nu_1(k-3)}
\]
and the  estimates for $\Psi_{jk}$ follow.

\noindent Similarly,
\begin{eqnarray*}\ab{\Psi_{jk}'}&\leqslant&\sum_{i=j}^k\ab{\wE(Y_j,\dots,Y_i',\dots,Y_k)}\leqslant
\sum_{i=j}^k2^{k-j}\sqrt{\Expect\ab{Y'_i}^2}\prod_{l\ne
i}^k\sqrt{\Expect\ab{Y_l}^2}\\
&\leqslant& 2^{k-j-1}(k-j+1)\prod_{i=j}^k\sqrt{a^2\nu_1(i)}
\end{eqnarray*}
and hence, the remaining two estimates follow. \qed
%
\begin{lemma}\label{fi} Let $a^2 \nu_1^{(n)} \leqslant 0.01$ and $s\geqslant 1$ be an integer. Then, for $k=1,2,\dots,n$ and  $t\in\RR$,
\begin{eqnarray}
\ab{\vfi_k(u)-1}&\leqslant&\frac{a^2}{6}[10\nu_1(k-1)+13\nu_1(k)],\quad
\ab{\vfi_k(u)-1}\leqslant \frac{1}{25},\quad \frac{1}{\ab{\vfi_k(u)}}\leqslant\frac{10}{9},\label{fi1}\\
\vfi_k(u)&=&1+\sum_{m=1}^s\frac{\nu_m(k)}{m!}(\ee^u-1)^m+\theta\bigg\{\frac{\ee^{hC_0}(\ee^h+1)^{s+1}
\nu_{s+1}(k)}{(s+1)!}\nonumber\\
&&+
(3.53)a^4\sum_{l=0}^3\nu_1^2(k-l)+(1.8)a^3[\nu_1^2(k-1)+\nu_1^2(k)]\nonumber\\
&&+(1.8)a^3\ab{\Cov(X_{k-1},X_k)}\bigg\},\label{fi2}\\
\ab{\vfi_k'(u)}&\leqslant&
2a^2[\nu_1(k-1)+\nu_1(k)],\quad\ab{\vfi_k'(u)}\leqslant
0.04,\label{fi3}\\
&&  \nonumber \\
\vfi_k'(u)&=&\ii \sum_{m=1}^s\frac{\nu_m(k)}{(m-1)!}(\ee^u-1)^{m-1}\ee^u+\theta\bigg\{\frac{\ee^{hC_0}(\ee^h+1)^{s}
\nu_{s+1}(k)}{s!}\nonumber\\
&&+(8.2)a^4\sum_{l=0}^3\nu_1^2(k-l)+2.6a^3[\nu_1^2(k-1)+\nu_1^2(k)]\nonumber\\
&&+ (2.6)a^3\ab{\Cov(X_{k-1},X_k)}\bigg\}.\label{fi4}
\end{eqnarray}
\end{lemma}

\noindent \Proof Further on, we assume that $k\geqslant 4$. For
smaller values of $k$, all the proofs indeed become shorter. For
brevity, we omit  the argument $u$, whenever possible. First note
that for all $t\in\RR$, $u\in\eK$. Indeed,
\[w(u)=\max_j\sqrt{\Expect\ab{Y_j}^2}\leqslant\max_j\sqrt{a^2\nu_1(j)}\leqslant
\frac{1}{10}.\]
 Consequently, by (\ref{fir}) and (\ref{Yj})
\begin{eqnarray*}\ab{\vfi_k-1}&\leqslant&\Expect\ab{Y_k}+\frac{2\Big(\Expect\ab{Y_{k-1}}^2\Expect\ab{Y_k}^2\Big)^{1/2}}{1-4w(u)} \nonumber \\
&\leqslant &
a\nu_1(k)+\frac{10}{3}(a^4\nu_1(k-1)\nu_1(k))^{1/2}\\
&\leqslant&
\frac{a^2\nu_1(k)}{2}+\frac{5a^2}{3}[\nu_1(k-1)+\nu_1(k)] \nonumber \\
 &=&
\frac{a^2}{6}[10\nu_1(k-1)+13\nu_1(k)].
\end{eqnarray*}
Using the assumption and noting that $1/\ab{\vfi_k}\leqslant
1/(1-\ab{\vfi_k-1})$, we obtain (\ref{fi1}). By (\ref{repr})
 \begin{equation}
\vfi_k=1+\Expect
Y_k+\frac{\Psi_{k-1,k}}{\vfi_{k-1}}+\frac{\Psi_{k-2,k}}{\vfi_{k-2}\vfi_{k-1}}+
\sum_{j=1}^{k-3}\frac{\Psi_{j,k}}{\vfi_j\cdots\vfi_{k-1}}.\label{skleid}
\end{equation}
Using Lemma \ref{Psi}, it follows that
\begin{equation}\label{o1}
\sum_{j=1}^{k-3}\frac{\ab{\Psi_{jk}}}{\vfi_j\dots\vfi_{k-1}}\leqslant
250a^4\sum_{l=0}^3\nu_1^2(k-l)\sum_{j=1}^{k-3}\bigg(\frac{10}{9}\bigg)^{k-j}
\bigg(\frac{1}{5}\bigg)^{k-j}\leqslant
(3.53)a^4\sum_{l=0}^3\nu_1^2(k-l).\end{equation}
 Similarly,  we have from (\ref{Yj})
 \begin{eqnarray*}
 \ab{\wE(Y_{k-1},Y_k)}&\leqslant&\Expect\ab{Y_{k-1},Y_k}+\Expect\ab{Y_{k-1}}\Expect\ab{Y_k}\leqslant
 a^2\Expect
 X_{k-1}X_k+a^2\nu_1(k-1)\nu_1(k)\\
 &=&a^2\Cov(X_{k-1},X_k)+a^22\nu_1(k-1)\nu_1(k)\\
&\leqslant&
 a^2\ab{\Cov(X_{k-1},X_k)}+a^2[\nu_1^2(k-1)+\nu_1^2(k)].
 \end{eqnarray*}
 Due to the trivial estimate $a\geqslant 2$,
 \begin{equation}\frac{\ab{\Psi_{k-1,k}}}{\ab{\vfi_{k-1}}}\leqslant\frac{5a^3}{9}\Big(
\ab{\Cov(X_{k-1},X_k)}+[\nu_1^2(k-1)+\nu_1^2(k)]\Big).
\label{o2}\end{equation}
 By assumption,  $\nu_1^{(n)}\leqslant 1/400$ and
\begin{eqnarray*}
\ab{\wE(Y_{k-2},Y_{k-1},Y_k)}&\leqslant&
[\ee^{hC_0}(2+h)]^3\big(\Expect
(X_{k-2}X_{k-1}X_k)+\nu_1(k-2)\Expect (X_{k-1}X_k)\\
&&+\Expect
(X_{k-2}X_{k-1})\nu_1(k)+\nu_1(k-2)\nu_1(k-1)\nu_1(k)\big)\\
&\leqslant&[\ee^{hC_0}(2+h)]^3\big(\Expect
(X_{k-1}X_k)(C_0+1/400)+\nu_1(k-1)\nu_1(k)(C_0+1/400)\big)\\
&\leqslant&\frac{401a^3}{400}(\Expect
(X_{k-1}X_k)+\nu_1(k-1)\nu_1(k))\\
&\leqslant&
\frac{401a^3}{400}(\ab{\Cov(X_{k-1},X_k)}+[\nu_1^2(k-1)+\nu_1^2(k)]).\end{eqnarray*}
Therefore,
\begin{equation}\frac{\ab{\Psi_{k-2,k}}}{\ab{\vfi_{k-2}\vfi_{k-1}}}\leqslant\bigg(\frac{10}{9}\bigg)^2
\frac{401a^3}{400} \Big(\ab{\Cov(X_{k-1},X_k)}+[\nu_1^2(k-1)+\nu_1^2(k)]\Big).
\label{o3}\end{equation}
The proof of
(\ref{fi2}) now follows by combining the last estimate with
(\ref{o1}), (\ref{o2}), (\ref{skleid}) and (\ref{Yj1}).

\noindent We prove
 (\ref{fi3}) by induction. We have
\[\vfi'_k=\Expect
Y_k'+\sum_{j=1}^{k-1}\frac{\Psi_{jk}'}{\vfi_j\cdots\vfi_k}-\sum_{j=1}^{k-1}\frac{\Psi_{jk}}{\vfi_j\cdots\vfi_k}
\sum_{m=j}^{k-1}\frac{\vfi_m'}{\vfi_m}.\]
 Applying Lemma \ref{Psi} and using (\ref{Yj}) and (\ref{Yj0}), we then get
 \begin{eqnarray*}\ab{\vfi_k'}&\leqslant&
 \frac{a^2}{4}\nu_1(k)+\frac{5}{2}a^2[\nu_1(k-1)+\nu_1(k)]\sum_{j=1}^{k-1}\bigg(\frac{2}{9}\bigg)^{k-j}(k-j+1)\\
 &&+
 5a^2[\nu_1(k-1)+\nu_1(k)]\frac{10}{9}(0.04)\sum_{j=1}^{k-1}(k-j)\bigg(\frac{2}{9}\bigg)^{k-j}\\
 &\leqslant&
 a^2[\nu_1(k-1)+\nu_1(k)]\bigg(\frac{1}{4}+\frac{80}{49}+\frac{4}{49}\bigg)\\
 & \leqslant &
 2a^2[\nu_1(k-1)+\nu_1(k)].
 \end{eqnarray*}
The proof of (\ref{fi4}) is similar to the proof of (\ref{fi2}).
We have \begin{eqnarray}\ab{\vfi'_k-\Expect
Y_k'}&\leqslant&\sum_{j=1}^{k-3}\bigg(\frac{10}{9}\bigg)^{k-j}\ab{\Psi_{jk}'}
+\sum_{j=1}^{k-3}\bigg(\frac{10}{9}\bigg)^{k-j}\ab{\Psi_{jk}}(k-j)\bigg(\frac{2}{45}\bigg)\nonumber\\
&&
+\sum_{j=k-2}^{k-1}\bigg(\frac{10}{9}\bigg)^{k-j}\ab{\Psi_{jk}'}
+\sum_{j=k-2}^{k-1}\frac{\ab{\Psi_{jk}}}{\ab{\vfi_j\cdots\vfi_{k-1}}}(k-j)\bigg(\frac{2}{45}\bigg).
\label{isv1}
\end{eqnarray}
Applying Lemma \ref{Psi}, we prove that
 \begin{equation}\label{isv2}
\sum_{j=1}^{k-3}\bigg(\frac{10}{9}\bigg)^{k-j}\ab{\Psi_{jk}'}
+\sum_{j=1}^{k-3}\bigg(\frac{10}{9}\bigg)^{k-j}\ab{\Psi_{jk}}(k-j)\bigg(\frac{2}{45}\bigg)
\leqslant (8.2)a^4\sum_{l=0}^3\nu_1^2(k-l).\end{equation}
 From (\ref{o2}) and (\ref{o3}), it follows that
 \begin{equation}
\sum_{j=k-2}^{k-1}\frac{\ab{\Psi_{jk}}}{\ab{\vfi_j\cdots\vfi_{k-1}}}(k-j)\bigg(\frac{2}{45}\bigg)\leqslant
(0.135)a^3(\ab{\Cov(X_{k-1},X_k)}+[\nu_1^2(k-1)+\nu_1^2(k)]).\label{isv3}
 \end{equation}
Taking into account (\ref{Yj}) and (\ref{Yj0}), we obtain
\[\ab{\wE(Y_{k-1}',Y_k)}\leqslant\Expect\ab{Y_{k-1}'Y_k}+\Expect\ab{Y_{k-1}'}
\Expect\ab{Y_k}\leqslant\ee^{2hC_0}(2+h)\Expect
(X_{k-1}X_k)+\frac{a^2}{2}\nu_1(k-1)\nu_1(k)\] and
\begin{equation}
\frac{10\ab{\Psi_{k-1,k}'}}{9}\leqslant
\frac{10a^2}{9}\Big(\ab{\Cov(X_{k-1},X_k)}+[\nu_1(k-1)+\nu_1(k)]\Big).
\label{isv4}
\end{equation}
Similarly,
\[
\bigg(\frac{10}{9}\bigg)^{2}\ab{\Psi_{k-2,k}'}\leqslant
(1.86)a^3(\ab{\Cov(X_{k-1},X_k)}+[\nu_1(k-1)+\nu_1(k)]).
\]
Combining the last estimate with (\ref{isv1})-(\ref{isv4}) and
(\ref{Yj2}), we complete the proof of (\ref{fi4}). \qed

\vspace*{0.4cm}
Let now
 \begin{equation} \label{neqn34}
 A(u)=\sum_{k=1}^n\ln\vfi_k(u)=\sum_{k=1}^n\sum_{j=1}^\infty\frac{(-1)^{j+1}(\vfi_k(u)-1)^j}{j}.
 \end{equation}
\begin{lemma}\label{A} Let $a^2\nu_1^{(n)}\leqslant 1/100$. Then
for all $t\in\RR$,
\begin{eqnarray}
\ab{A}&\leqslant&4a^2\Gamma_1,\quad\ab{A'}\leqslant 4a^2\Gamma_1,\label{A1}\\
A&=&
\sum_{m=1}^s\frac{(\ee^u-1)^m}{m!}\sum_{k=1}^n\nu_m(k)+\theta\bigg\{\frac{\ee^{hC_0}(\ee^h+1)^{s+1}}{(s+1)!}
\sum_{k=1}^n\nu_{s+1}(k)\nonumber\\
&&+24a^4\sum_{k=1}^n\nu_1^2(k)+(1.8)a^3\sum_{k=2}^n\ab{\Cov(X_{k-1},X_k)}\bigg\},\label{A2}\\
A'&=&\ii\sum_{m=1}^s\frac{\nu_m(k)}{(m-1)!}(\ee^u-1)^{m-1}\ee^u+\theta\bigg\{\frac{\ee^{hC_0}(\ee^h+1)^{s}
}{s!}\sum_{k=1}^n\nu_{s+1}(k)\nonumber\\
&&+(51.4)a^4\sum_{k=1}^n\nu_1^2(k)+
(2.6)a^3\sum_{k=2}^n\ab{\Cov(X_{k-1},X_k)}\bigg\}.\label{A3}
\end{eqnarray}
\end{lemma}

\noindent \Proof Using the first estimate in (\ref{fi1}), we have
$\ab{\vfi_k-1}\leqslant 0.04$. Therefore,
\begin{equation*}
\ab{A}\leqslant\sum_{k=1}^n\ab{\vfi_k-1}\sum_{j=1}^\infty
(0.04)^{j-1} \leqslant
\Big(\frac{1}{0.96}\Big)\sum_{k=1}^n\frac{a^2[10\nu_1(k-1)+13\nu_1(k)]}{6}\leqslant
4a^2\Gamma_1.\label{pr2}
\end{equation*}
Similarly,
 \[\ab{A'}\leqslant\sum_{k=1}^n\frac{\ab{\vfi_k'}}{\ab{\vfi_k}}\leqslant\frac{10}{9}\sum_{k=1}^n\ab{\vfi_k'}\leqslant
\frac{20a^2}{9}\sum_{k=1}^n[\nu_1(k-1)+\nu_1(k)]\leqslant
4a^2\Gamma_1.\]
 From Lemma \ref{fi}, it follows
\[
\ab{\vfi_k-1}^2\leqslant\frac{a^4}{36}\Big(10\nu_1(k-1)+13\nu_1(k)\Big)^2\leqslant\frac{a^4}{36}\Big(230\nu_1^2(k-1)+
299\nu_1^2(k)\Big)\] and
\[
\sum_{k=1}^n\sum_{j=2}^\infty\frac{\ab{\vfi_k-1}^{j-2}}{j}\leqslant\frac{1}{2}\sum_{k=1}^n\ab{\vfi_k-1}^2\sum_{j=2}^\infty
(0.04)^{j-2}\leqslant (7.66)a^4\sum_{k=1}^n\nu_1^2(k).
\]
Consequently,
\[
A=\sum_{k=1}^n(\vfi_k-1)+ (7.66) \theta a^4\sum_{k=1}^n\nu_1^2(k)
\]
and (\ref{A2}) follows from Lemma \ref{fi} and the rough estimate
$a^3\leqslant a^4/2$, since $a\geqslant 2$.

\noindent For the proof of (\ref{A3}), note that
\[A'=\sum_{k=1}^n\vfi_k'+\sum_{k=1}^n\frac{\vfi_k'}{\vfi_k}(1-\vfi_k)\]
and applying a slightly sharper estimate than in Lemma \ref{fi},
namely $1/\ab{\vfi}\leqslant 25/24$, we obtain
\[\sum_{k=1}^n\frac{\ab{\vfi_k'}}{\ab{\vfi_k}}\ab{1-\vfi_k}\leqslant
\frac{25a^4}{72}\sum_{k=1}^n[\nu_1(k-1)+\nu_1(k)][10\nu_1(k-1)+13\nu_1(k)]\leqslant
16a^4\sum_{k=1}^n\nu_1^2(k).
\]
Now, it remains to apply (\ref{fi4}) to complete the proof. \qed


\begin{lemma} \label{apvertimas}  Let $M$ be a  finite variation measure concentrated on
integers and  $\sum_k\ab{k}\ab{M\{k\}}<\infty$. Then,
\begin{equation*}
 \sum_{k=-\infty}^\infty\ab{M\{k\}}\leqslant\bigg(\frac{1}{2}+\frac{1}{2\pi}\bigg)^{1/2}\biggl(
\int\limits_{-\pi}^{\pi}\ab{\w M(\ii t)}^2+ \ab{\w M'(\ii t)}^2
\dd t\biggr)^{1/2}. \label{TVAP}
 \end{equation*}
\end{lemma}
Lemma \ref{apvertimas} is  a special case of the Lemma  from
Presman (1985).

In the following lemma,  we assume that
$M(k)=M\{[0,k]\}=\sum_{j=0}^kM\{j\}$.
\begin{lemma} \label{WSJ} Let  $g(k)>0$ be an increasing function
 and $M$ be a measure,  both defined on nonnegative
integers. Assume $\sum_{k=0}^\infty M\{k\}=0$ and
$\sum_{k=0}^\infty kg(k)\ab{M\{k\}}<\infty$. Then,
\[
\sum_{k=0}^\infty g(k)\ab{M(k)}\leqslant\sum_{k=0}^\infty
kg(k)\ab{M\{k\}}.
\]
If $g(k)=\ee^{hk}$, for some  $h>0$, then
\[
\sum_{k=0}^\infty
\ee^{hk}\ab{M(k)}\leqslant\frac{1}{\ee^h-1}\sum_{k=0}^\infty\ee^{hk}\ab{M\{k\}}.
\]
\end{lemma}

\noindent \Proof We have
\begin{eqnarray*}
\sum_{k=0}^\infty g(k)\ab{M(k)}&=&\sum_{k=0}^\infty
g(k)\Ab{\sum_{j=0}^kM\{j\}}= \sum_{k=0}^\infty
g(k)\Ab{\sum_{j=k+1}^\infty M\{j\}}
\leqslant\sum_{k=0}^\infty g(k)\sum_{j=k+1}^\infty\ab{M\{j\}}\\
&=&\sum_{j=1}^\infty\ab{M\{j\}}\sum_{k=0}^{j-1} g(k)\leqslant
\sum_{j=0}^\infty jg(j)\ab{M\{j\}}.
\end{eqnarray*}
When $g(k)=\ee^{hk}$,  $h>0$, we simply note that
\begin{equation*}
\sum_{k=0}^{j-1}\ee^{hk}=\frac{\ee^{hj}-1}{\ee^h-1}<\frac{\ee^{hj}}{\ee^h-1}.
~~~~~~~~~~~~~~~~~~~~~~~~~~~~~~~~~~~~~~~~~~~~~~~~~~~~\rlap{$\qquad \Box$}
\end{equation*}

\noindent \textbf{Proof of Theorem \ref{teorema}.} Let  $M$
 defined by

\vspace*{-0.4cm}
\begin{equation*}
M\{k\}=\ee^{hk}(F_n\{k\}-\CP(s,\lambda)\{k\})
\end{equation*}

\vspace*{-0.4cm} be a (signed) measure. Then,
\[\w M(\ii t)=\sum_{k=0}^\infty
(F_n\{k\}-\CP(s,\lambda)\{k\})\ee^{hk+\ii tk}=\w F_n(\ii
t+h)-\w\CP(s,\lambda)(\ii t+h).\]
  Applying Lemma \ref{apvertimas},  we obtain
\begin{equation}
\sum_{k=0}^\infty\ee^{hk}\ab{F_n\{k\}-\CP(s,\lambda)\{k\}}\leqslant\sqrt{\pi+1}\,\bigl(\sup_{\ab{t}\leqslant\pi}\ab{\w
M(\ii t)}+\sup_{\ab{t}\leqslant\pi}\ab{\w {M^{'}}(\ii t)}\bigr).
\label{pr1}
\end{equation}

\noindent Let $A$ be defined as in  (\ref{neqn34}) and
\vspace*{-0.3cm}
\[
B=\sum_{j=1}^s\lambda_j(\ee^{ju}-1)=\sum_{m=1}^s(\ee^u-1)^m\sum_{j=m}^s\binomial{j}{m}\lambda_j.\]
Then, $\ab{B}\leqslant\sum_{j=1}^s\lambda_j(\ee^{jh}+1).$ Also,
$\w M(\ii t)=\ee^A-\ee^B$ and
 \begin{eqnarray}
 \ab{\w M(\ii
t)}&\leqslant&\ab{\ee^A-\ee^B}\leqslant\psi\ab{A-B}, \label{pro1}\\
\ab{\w M'(\ii
t)}&\leqslant&\ab{A'}\ab{\ee^A-\ee^B}+\ab{\ee^B}\ab{A'-B'}\leqslant
\psi\ab{A'}\ab{A-B}+\psi\ab{A'-B'},\label{pro2}
\end{eqnarray}
where $\psi$ is defined in (\ref{eqn11}).

\noindent The inequality $\ab{\ee^A-\ee^B}\leqslant\psi\ab{A-B}$ follows
from the fact that if the real part of a complex number $Re
z\leqslant 0$, then
\[\ab{\ee^z-1}\leqslant\Ab{\int_0^1 z\ee^{\tau
z}\dd\tau}\leqslant\ab{z}\int_0^1\exponent{\tau Re
z}\dd\tau\leqslant\ab{z}.\]
 Indeed, if $Re(A-B)<0$, then
 \[\ab{\ee^A-\ee^B}=\ab{\ee^B}\ab{\ee^{A-B}-1}\leqslant\ab{\ee^B}\ab{A-B}\leqslant\psi\ab{A-B}.\]
 If $Re(B-A)\leqslant 0$, then
 \[\ab{\ee^A-\ee^B}=\ab{\ee^A}\ab{1-\ee^{B-A}}\leqslant\ab{\ee^A}\ab{A-B}\leqslant\psi\ab{A-B}.\]

\noindent The proof is now completed by combining (\ref{pro1}),
(\ref{pro2}) with (\ref{pr1}) and using Lemma \ref{A}. \hfill{\qed}

\subsection*{Acknowledgment}
The authors wish to thank the referees for  helpful comments which
helped to improve the paper.



\end{document}